\documentclass[11pt]{article}

\usepackage{amssymb}
\usepackage{amsfonts,amstext,amsmath,amsthm,
latexsym,mathrsfs,amsbsy,mparhack}
\usepackage[bbgreekl]{mathbbol}

\usepackage{graphicx}
\usepackage{color}
\usepackage{bold-extra}

\usepackage{marginnote}

\usepackage[russian,english,greek]{babel}

\usepackage[b]{esvect}

\usepackage{tocloft}

\usepackage{doi}

\hypersetup
{%
hypertex=true,
hypertexnames=true,
hidelinks,
pagebackref=true,
}

\usepackage
{lineno}

\input{69.sty}

\begin{document}

\selectlanguage{english}


\title{The parameterfree Comprehension does not imply the 
full Comprehension in the 2nd order Peano arithmetic\thanks
{This paper was written under the support 
of RFBR (Grant no~20-01-00670).}
}

\author{Vladimir Kanovei\thanks
{Institute for Information Transmission Problems
(Kharkevich Institute) of Russian Academy
of Sciences (IITP), Moscow, Russia, 
{\tt kanovei@iitp.ru} } 
\and 
Vassily~Lyubetsky\thanks
{Institute for Information Transmission Problems
(Kharkevich Institute) of Russian Academy
of Sciences (IITP), Moscow, Russia, 
{\tt lyubetsk@iitp.ru} } 
}

\date{\today}
\maketitle


\begin{abstract}
The parameter-free part $\pads$ 
of $\pad$, the 2nd order Peano arithmetic, is 
considered.
We make use of a product/iterated Sacks forcing  
to define an $\omega$-model of $\pads + \CA(\us12)$, 
in which an example of the full Comprehension 
schema $\CA$ fails.
Using Cohen's forcing, we also define an $\omega$-model 
of $\pads$, in which not every set has its complement, 
and hence the full $\CA$ fails in a rather elementary way.
\end{abstract}

\vyk{
}


\setcounter{tocdepth}{1}

\setlength{\cftbeforesecskip}{5.5pt}

{\footnotesize\def\contentsname{\normalsize \bf
Contents
}\tableofcontents}



\np


\parf{Introduction}
\las{int}

Discussing the structure and deductive properties of 
the second order Peano arithmetic $\pad$, Kreisel 
\cite[\S\,III, page 366]{Kr} wrote that 
the selection of subsystems ``is a central problem''. 
In particular, Kreisel notes, that 
\begin{quote}
[...] if one is convinced of the significance 
of something like a given axiom schema, 
it is natural to study details, such as the effect 
of parameters. 
\end{quote}
Recall that 
\rit{parameters} in this context are free variables 
in various axiom schemata in $\PA$, $\ZFC$, 
and other similar theories. 
Thus the most obvious way to study 
``the effect of parameters'' is to compare the strength 
of a given axiom schema $S$ with its parameter-free 
subschema $S^\ast$. 
(The asterisk will mean the parameter-free subschema in this 
paper.) 

Some work in this direction was done in the early years of 
modern set theory. 
In particular Guzicki~\cite{guz} proved that the Levy-style 
generic collapse 
(see, \eg, Levy~\cite{levy2} and Solovay~\cite{sol})
of all cardinals $\om_\al^\rL$, 
$\al<\om_1^\rL$, results in a generic extension of $\rL$ in 
which the (countable) choice schema $\AC$, 
in the language of $\pad$, fails but its \paf\ subschema 
$\ACs$ holds, so that $\ACs$ is strictly weaker than $\AC$. 
This can be compared with an opposite result for the 
\rit{dependent choice} schema $\DC$, 
in the language of $\pad$, which is equivalent to 
its \paf\ subschema $\DCs$ by a simple argument given in 
\cite{guz}. 

Some results related to \paf\ versions of the Separation 
and Replacement axiom schemata in $\zfc$ 
also are known from \cite{corrada,levyp,SS}. 

This paper is devoted to the role of parameters in the 
\rit{comprehension schema} $\CA$ of $\pad$. 
Let $\pads$ be the subtheory of\/ $\pad$ in 
which the full schema $\CA$ is replaced by its 
\paf\ version $\CAs$, 
and the Induction principle is formulated as 
a schema rather than one sentence. 
The following 
Theorems \ref{mt1} and \ref{mt2} are our main results. 

\bte
\lam{mt1}
Let\/ $\coh$ be the Cohen forcing for adding a generic 
subset of\/ $\omega$. 
Let\/ $\cohw$ be the finite-support product.
Suppose that\/ $\sis{x_i}{i<\omi}$ is a sequence\/ 
\dd\cohw generic over\/ $\rL$, the constructible 
universe. 

Let\/ $X=(\pws\om\cap\rL)\cup \ens{x_i}{i<\om}$.
Then\/ $\stk{\om}{X}$ is a model 
of\/ $\pads$, 
but not a model of\/ $\CA$ as\/ $X$ does not 
contain the complements\/ $\om\bez x_i$. 

Thus\/ $\CA$, even in the particular form 
claiming that every set 
has its complement, is not provable in\/ $\pads$.
\ete

It is quite obvious that a subtheory like $\pads$, 
that does not allow such a fundamental thing as the 
complement formation, is unacceptable. 
This is why we adjoin $\CA(\us12)$, \ie, the full $\CA$ 
(with parameters) restricted to $\us12$ formulas, 
in the next theorem, 
to obtain a more plausible subsystem. 

\bte
\lam{mt2}
There is a generic extension\/ $\rL[G]$ of\/ $\rL$ 
and a set\/ $M\in\rL[G]$, such that\/ 
$\pws\om\cap\rL\sq M\sq\pws\om$ and\/ 
$\stk{\om}{M}$ is a model of\/ $\pads + \CA(\us12)$ 
but not a model of\/ $\pad$.

Therefore\/ $\CA$ is not provable even in\/ 
$\pads + \CA(\us12)$.
\ete

Theorem~\ref{mt2} will be established by means of 
a complex product/iteration of the Sacks forcing 
and the associated coding by degrees of 
constructibility, approximately as discussed in 
\cite[page~143]{matsur}, around Theorem~T3106.

Identifying the theories with their deductive 
closures, we may present the concluding statements 
of Theorems \ref{mt1} and \ref{mt2} as resp.
\bce
$\pads\sneq\pad$ \ 
and \ $(\pads + \CA(\us12))\sneq\pad$.
\ece
Studies on subsystems of $\pad$ have discovered many 
cases in which $S\sneq S'$ holds for a given pair 
of subsystems $S,S'$, see \eg\ \cite{simp}. 
And it is a rather typical case that such a strict 
extension is established by demonstrating that $S'$ 
proves the consistency of $S$.
One may ask whether this is the case for the 
results in the displayed line above. 
The answer is in the negative: 
namely the theories $\pads$, $\pads + \CA(\us12)$, 
and the full $\pad$ happen to be   
equiconsistent by a result in \cite{HFuse81}, 
also mentioned in \cite{schindt}. 
This equiconsistency result also follows from 
a somewhat sharper 
theorem in \cite[1.5]{Schm}.\snos
{The authors are thankful to Ali Enayat for the 
references to \cite{HFuse81,schindt,Schm} in 
matters of this equiconsistency result.}

\parf{Preliminaries} 
\las{prel}

Following \cite{aptm,Kr,simp} we define the second order 
Peano arithmetic $\pad$ as a theory in the language 
$\lpad$ with two sorts of variables -- 
for natural numbers and for sets of them. 
We use $j,k,m,n$ for variables over $\omega$ and 
$x,y,z$ for variables over $\pws\om$, reserving capital 
letters for subsets of $\pws\om$ and other sets. 
The axioms are as follows:
\ben
\nenu
\itlb{pa1}%
Peano's axioms for numbers.

\itlb{pa2}%
The Induction schema 
$\Phi(0) \land \kaz k\,(\Phi(k)\imp\Phi(k+1))
\limp \kaz k\,\Phi(k)$, 
for every formula $\Phi(k)$ in $\lpad$, 
and in $\Phi(k)$ we allow parameters, 
\ie, free variables other than $k$.\snos
{We cannot use Induction as one sentence because 
the Comprehension schema $\CA$ is not  
assumed in full generality in the context of 
Theorem~\ref{mt1}.}

\itlb{pa3}%
Extensionality for sets.
  
\itlb{pa4}%
The Comprehension schema $\CA$:
$\sus x \,\kaz k\,(k\in x\eqv\Phi(k))$, 
for every formula $\Phi$ in which the variable $x$ 
does not occur, and in $\Phi$ we allow parameters.
\een

We let $\CA(\us12)$ be the full $\CA$ 
restricted to $\us12$ formulas $\Phi$.\snos
{A $\us12$ formula is any $\lpad$ formula of the form 
$\kaz x\,\sus y\,\Psi$, where $\Psi$ does not contain 
quantified variables over $\pws\om$.} 

We let $\CAs$ be the \paf\ sub-schema of $\PA$ 
(that is, $\Phi(k)$ contains no free variables 
other than $k$). 

We let $\pads$ be the subsistem of $\pad$ with 
$\CA$ replaced by $\CAs$. 

\bre
\lam{r1}
In spite of Theorem~\ref{mt1}, $\pads$ proves 
$\CA$ with parameters over $\om$ 
(but not over $\pws\om$) allowed. 
Indeed suppose that $\Phi$ is $\Phi(k,m)$ in 
\ref{pa4} and $\Phi$ has no other free variables. 
Arguing in $\pad$, assume towards the contrary 
that the formula 
$\psi(m):=\sus x \,\kaz k\,(k\in x\eqv\Phi(k,m))$
holds not for all $m$. 
By Induction, take the least $m$ for which 
$\psi(m)$ \rit{fails}. 
This $m$ is definable, and therefore it can be 
eliminated, and hence we have $\psi(m)$ for this 
$m$ by $\CAs$. 
This is a contradiction.
\ere


\parf{Extension by Cohen reals} 
\las{coh}

Here we prove Theorem~\ref{mt1}. 
We assume some knowledge 
of forcing and generic models, as \eg\ in 
Kunen~\cite{kun}, especially Section IV.6 there 
on the ``forcing over the universe'' approach.

Recal that the Cohen forcing notion $\coh=\bse$ 
consists of all finite dyadic tuples including 
the empty tuple $\La$. 
If $u,v\in\bse$ then $u\su v$ means that $v$ is 
a proper extension of $u$, whereas $u\sq v$ means 
$u\su v\lor u=v$.
The \rit{finite-support product} $\bP=\bsep$ consists 
of all maps $p:\om\to\bse$ such that $p(i)=\La$ 
(the empty tuple) for all but finite $i<\om$. 
The set $\bP$ is ordered opposite to the componentwise 
extension, so that $p\leq q$ 
($p$ is \rit{stronger} as a forcing condition) 
iff $q(i)\sq p(i)$ for all $i<\om$.
The condition $\lao$ defined by $\lao(i)=\La$, 
$\kaz i$, is the \dd\leq largest 
(the weakest) 
element of $\bP$.

We consider the set $\Perm$ of all 
idempotent 
\rit{permutations} of $\om$, that is, 
all bijections $\pi:\om\onto\om$ such that 
$\pi=\pi\obr$ and 
\rit{the domain of nontriviality} 
$\abc \pi=\ens{i}{\pi(i)\ne i}$ is finite. 
If $\pi\in\Perm$ and $p$ is a function with 
$\dom\pi=\om$, then $\pi p$ is defined by 
$\dom (\pi p)=\om$ and $(\pi p)(\pi(i))=p(i)$
for all $i<\om$, so formally 
$\pi p=p\circ\pi\obr=p\circ \pi$ 
(the superposition). 
In particular if $p\in\bP$ then $\pi p\in\bP$ 
and 
$\abc{\pi p}=\ima \pi {\abs p}=
\ens{\pi(i)}{i\in\abc p}$.

\bpf[Theorem~\ref{mt1}]
We make use of G\"odel's \rit{constructible universe} 
$\rL$ as \rit{the ground model} for our forcing 
constructions. 
Suppose that $G\sq\bP$ is a set \dd\bP generic over 
$\rL$. 
If $i<\om$ then 
\bit
\item[$-$]\msur 
$G_i=\ens{p(i)}{p\in G}\sq\bse$ 
is a set \dd\bse generic (Cohen generic) over~$\rL$, 

\item[$-$]\msur 
$a_i[G]=\bigcup G_i\in\dn$ is a real 
Cohen generic over $\rL$, and 

\item[$-$]\msur 
$x_i[G]=\ens{n}{a_i(n)=1}\sq\om$ 
is a subset of $\om$ Cohen generic over $\rL$. 

\item[$-$]\msur 
$X=X[G]=(\pws\om\cap\rL)\cup\ens{x_i[G]}{i<\om}$. 
\eit
Thus $X[G]\in\rL[G]$ and $X[G]$ consists of all 
subsets of $\om$ already in $\rL$ 
and all Cohen-generic sets $x_i[G]$, $i<\om$. 

We assert that the model $\stk\om\xg$ 
proves Theorem~\ref{mt1}. 

The only thing to check is that $\stk\om\xg$ 
satisfies $\CAs$. 
For that purpose, assume that $\Phi(k)$ is a 
\paf\ $\lpad$ formula with $k$ the only 
free variable. 
Consider the set 
$y=\ens{k<\om}{\stk\om\xg\mo\Phi(k)}$; 
then $y\in\rL[G]$, $y\sq\om$. 
We claim that in fact $y$ 
belongs to $\rL$, and hence to $\xg$. 

Let $\ofo$ be the forcing relation associated 
with $\bP$. 
In particular, if $p\in\bP$ and $\psi$ 
is a parameter-free formula then 
$p\ofo\psi$ iff $\psi$ holds in any 
\dd\bP generic extension $\rL[H]$ of $\rL$ 
such that $p\in H$. 

Let $\uG$ be a canonical \dd\bP name for $G$.
We assert that 
\busq
{eq1}
{
y=\ens{k<\om}
{\lao\ofo\text{``$\stk\om\xug\mo \Phi(k)$''}}.
}
Indeed assume that the condition $\lao$ \dd\bP 
forces ``$\stk\om\xug\mo \Phi(k)$''. 
But $\lao\in G$ since $\lao$ is the weakest condition 
in $\bP$. 
Therefore 
$\stk\om\xg\mo \Phi(k)$ by the forcing theorem, 
thus $k\in y$, as required.

To prove the converse, assume that $k\in y$. 
Then by the forcing theorem there is a condition 
$p\in G$ forcing ``$\stk\om\xug\mo \Phi(k)$''. 
We claim that then $\lao$ forces the same as well. 

Indeed otherwise there is a condition $q\in\bP$ 
which forces ``$\stk\om\xug\mo \neg\,\Phi(k)$''. 
There is a permutation $\pi\in\Perm$ satisfying 
$\abc{r}\cap\abc p=\pu$, where $r=\pi q\in\bP$. 
We claim that $r$ forces 
``$\stk\om\xug\mo \neg\,\Phi(k)$''. 
Indeed assume that $H\sq\bP$ is a set \dd\bP generic 
over $\rL$, and $r\in H$. 
We have to prove that 
$\stk\om{\xh}\mo\neg\,\Phi(k)$.
The set $K=\ens{\pi{r'}}{r'\in H}$ is \dd\bP generic 
over $\rL$ along with $H$ since $\pi\in\rL$. 
Moreover $K$ contains $q$. 
It follows that 
$\stk\om{\xk}\mo\neg\,\Phi(k)$
by the forcing theorem and the choice of $q$. 
However the sequence $\sis{x_i[K]}{i<\om}$ is 
equal to the permutation of 
the sequence $\sis{x_i[H]}{i<\om}$ by $\pi$. 
It follows that $\xh=\xk$, and hence 
$\stk\om{\xh}\mo\neg\,\Phi(k)$, as required. 
Thus indeed $r$ forces 
``$\stk\om\xug\mo \neg\,\Phi(k)$''. 

However $p$ forces ``$\stk\om\xug\mo \Phi(k)$'', 
and $p,r$ are compatible in $\bP$ because 
$\abc{r}\cap\abc p=\pu$. 
This is a contradiction. 

We conclude that $\lao$ forces 
$\stk\om\xug\mo \Phi(k)$, and this completes the 
proof of \eqref{eq1}. 

But it is known that the forcing relation $\ofo$ 
is expressible in $\rL$, the ground model. 
Therefore it follows from \eqref{eq1} 
that $y\in\rL$, hence $y\in \xg$, as required. 
\epf

\parf{Generalized Sacks iterations}
\las{gsi}

Here we begin the proof of Theorem~\ref{mt2}. 
The proof involves the engine of generalized 
product/iterated Sacks forcing developed in 
\cite{fm97,jsl99} on the base of earlier papers 
\cite{balav,groapp} and others. 
We still consider the constructible universe 
$\rL$ as the ground model for the extension, 
and define, in $\rL$, the set 
\busq
{eq2}
{\bI\;=\;(\omi\ti{\bse})\:\cup\:\omi\,;\quad 
\bI\in\rL\,,}
partially ordered so that 
$\ang{\ga,s}\cle\ang{\ba,t}$ iff $\ga=\ba$ and 
$s\sq t$ in $\bse,$ while the ordinals in $\omi$ 
(the second part of $\bI$) 
remain \dd\cle incomparable. 

Our plan is to define a  product/iterated 
generic Sacks extension $\rL[\va]$ of $\rL$ 
by an array $\va=\sis{a_\i}{\i\in\bI}$ of reals 
$a_\i\in\dn$, in which the 
structure of ``sacksness'' is determined by this 
set $\bI$, so that in particular each $a_\i$ is 
Sacks-generic 
over the submodel $\rL[\sis{a_\j}{\j\cls \i}]$. 

Then we define the set $\bJ\in\rL[\va]$ of all 
elements $\i\in\bI$ such that:\vtm 

--- either 
$\i=\ang{\ga,0^m}$, where $\ga<\omi$ and $m<\om$,\vom 

--- or $\i=\ang{\ga,0^m\we1}$, where $\ga<\omi$ and 
$m<\om$, $a_\ga(m)=1$.\vtm 

\noi
This any $\i=\ang{\ga,0^m}\in\bJ$ is a splitting 
node in $\bJ$ iff $a_\ga(m)=1$, or in other words
\busq
{eq3}
{a_\ga(m)=1\quad\text{iff}\quad
\ang{\ga,0^m}\text{ is a splitting node in }\bJ\,,}
We'll finally prove that the according set 
\busq
{eq4}
{M=\pws\om\cap\bigcup_{\i_1,\dots,\i_n\in\bJ} 
\rL[a_{\i_1},\dots,a_{\i_n}]
}   
leads to the model $\stk\om M$
for Theorem~\ref{mt2}.
The reals $a_\ga$ will not belong to $M$ 
by the choice of $\bJ$, but will be definable in 
$\stk\om M$ 
(with $a_{\ang{\ga,\La}}\sq\om$ as a parameter) 
via the characterization of the 
splitting nodes in $\bJ$ by \eqref{eq3}.


\parf{Iterated perfect sets}
\las{prelim}

{\ubf Arguing in $\rL$ in this section}, 
we define $\bI=\stk{\bI}{\cle}$ as above. 

Let $\cpo$ be the set of all countable 
(including finite) 
sets $\za\sq\bI$. 

If $\za\in\cpo$ then $\IS_\za$ is the set of all 
initial segments of $\za$. 

Greek letters $\xi,\,\eta,\,\za,\,\vt$ will 
denote sets in $\cpo$.  

Characters $\i,\,\j$ are used to denote 
{\it elements\/} of $\bI$.

For any $\i\in\za\in\cpo,$ 
we consider initial segments 
$\za\ile\i=\ans{\j\in\za:\j\cls\i}$ and 
$\za\nlq\i=\ans{\j\in\za:\j\not\cge\i},$  and 
$\za\ilq\i,\msur$ $\za\nle\i$ defined analogously. 

Further, $\bn$ is the {\it Baire space\/}. 
Points of $\bn$ will be called {\it reals\/}. 

Let $\can{}=\dn\sq\bn$ be the {\it Cantor space\/}. 
For any countable set 
$\xi,$ $\can\xi$ is the product of \dd\xi many 
copies of $\can{}$ with the product topology. 
Then every $\can\xi$ is a compact space, 
homeomorphic to $\can{}$ itself unless $\xi=\emps$.

Assume that $\eta\sq\xi\in\cpo$.  
If $x\in\can\xi$ then let 
$x\res\eta\in\can\eta$ denote the usual restriction. 
If $X\sq\can\xi$ then let 
$X\res \eta=\ans{x\res\eta:x\in X}$.
To save space, let $X\rsd{\cls\i}$ mean $X\res\xi\ile\i$, 
$X\rsd{\not\cge\i}$ mean $X\res\xi\nlq\i$, \etc 

But if $Y\sq\can\eta$ then we put 
$Y\ares \xi=\ans{x\in\can\xi:x\res\eta\in Y}$.

To describe the idea behind the definition of iterated 
perfect sets, recall that the Sacks forcing consists 
of perfect subsets of $\can{}$, that is, 
sets of the form 
$\ima H\cam =\ens{H(a)}{a\in\cam}$, 
where $H:\can{}\onto X$ is a homeomorphism. 

To get a product Sacks model, with two factors 
(the case of a two-element unordered set 
as the length of iteration), 
we have to consider sets $X\sq\can 2$ 
of the form $X= \ima H{\can 2}$ where $H$, 
a homeomorphism defined on 
$\can 2,$ splits in obvious way into a pair of 
one-dimentional homeomorphisms. 

To get an iterated Sacks model, 
with two stages of iteration  
(the case of a two-element ordered set 
as the length of iteration), 
we have to consider sets $X\sq\can 2$ 
of the form $X=\ima H{\can 2}$, 
where $H$, a homeomorphism defined on $\can 2,$ 
satisfies the following: if $H(a_1,a_2)=\ang{x_1,x_2}$ 
and $H(a'_1,a'_2)=\ang{x'_1,x'_2}$ then 
$a_1=a'_1\eqv x_1=x'_1$. 

The combined product/iteration case results 
in the following definition. 

\bdf
[iterated perfect sets, \cite{fm97,jsl99}]
\lam{ips}
For any 
$\za\in\cpo,$ $\pe\za$ is the collection of all 
sets $X\sq\can\za$ such that there is a homeomorphism 
$H:\can\za\onto X$ satisfying
\dm
x_0\res\xi=x_1\res\xi\,\eqv\,H(x_0)\res\xi=H(x_1)\res\xi
\dm
for all $x_0,\,x_1\in\dom H$ and $\xi\in\cpo$, 
$\xi\sq\za$. 
Homeomorphisms $H$ satisfying this requirement 
will be called {\it\prok\/}. 
In other words, sets in $\pe\za$ are images 
of $\can\za$ via \prok\ homeomorphisms. 
\edf 

\bre
\lam{lao}
Note that $\pu$, the empty set, formally belongs 
to $\cpo$, and then $\can\pu=\ans\pu$, and we 
easily see that $\bon=\ans\pu$ is the only set 
in $\pe\pu$. 
\ere

For the convenience of the reader, we now present 
five lemmas on sets in $\pe\za$ established in 
\cite{fm97,jsl99}.

\ble
[Proposition 4 in \cite{fm97}]
\lam{oldf}
Let $\za\in\cpo$.
Every set\/ $X\in\pe\za$ is closed and satisfies 
the following properties$:$ 
\ben
\itemsep=1mm
\def\theenumi{{\hskip1pt{\rm P-\arabic{enumi}}\hskip1pt}}
\itlb{perf} 
If\/ $\i\in\za$ and\/ $z\in X\rsd{<\i}$ 
then\/ 
$D_{Xz}(\i)=\ens{x(\i)}{x\in X\land x\rsd{<\i}=z}$ 
is a perfect set in\/ $\can{}$. 

\itlb{oz} 
If\/ $\xi\in\IS_\za$, and a set\/ $X'\sq X$ 
is open in\/ $X$ 
(in the relative topology) 
then the projection\/ $X'\res\xi$ is 
open in $X\res\xi$.
In other words, the projection from $X$ to $X\res\xi$ 
is an open map. 

\itlb{indep} 
If\/ $\xi,\eta\in\IS_\za$, 
$x\in X\res\xi$, $y\in X\res\eta$, and\/ 
$x\res (\xi\cap \eta)=y\res (\xi\cap \eta),$ then\/ 
$x\cup y\in X\res (\xi\cup \eta)$. 
\een
\ele
\bpf[sketch] 
Clearly $\can\za$ satisfies \ref{perf}, \ref{oz}, 
\ref{indep}, 
and one easily shows that \prok\  
homeomorphisms preserve the requirements.
\epf

\ble
[Lemma 6 in \cite{fm97}] 
\lam{less}
If\/ $\za\in\cpo$, $X\in\pe\za$, 
$\xi\in\IS_\za$, then\/ 
$X\res\xi\in\pe\xi$.
\ele

\ble
[Lemma 8 in \cite{fm97}] 
\lam{clop}
If\/ $\za\in\cpo$, $X\in\pe\za$, 
a set\/ $X'\sq X$ is open in\/ $X$, 
and\/ $x_0\in X',$ then there is a set\/  
$X''\in\pe\za$, $X''\sq X',$ clopen in\/ $X$ 
and containing\/ $x_0$.
\ele

\ble
[Lemma 10 in \cite{fm97}] 
\lam{apro}
Suppose that\/ $\za\in\cpo$, $\eta\in\IS_\za$,  
$X\in\pe\za$, $Y\in\pe\eta,$ and\/ $Y\sq X\res \eta$. 
Then\/ $Z=X\cap (Y\ares \za)$ belongs to\/ $\pe\za$. 
\ele

\ble
[Lemma 10 in \cite{jsl99}] 
\lam{99}
Suppose that\/ $\za\in\cpo$, $\xi\sq\za$,  
$X\in\pe\xi$. 
Then\/ $X\ares \za$ belongs to\/ $\pe\za$. 
\ele

\parf{The forcing and the basic extension}
\las{re}

This section introduces the forcing notion we 
consider and the according generic extension 
called the basic extension.

{\ubf We continue to argue in $\rL$.}  
Recall that a partially ordered set $\bI\in \rL$ 
is defined by \eqref{eq2} in Section~\ref{gsi}, and 
$\cpo$ is the set of all at most countable initial 
segments $\xi\sq\bI$ in $\rL$.
For any $\za\in\ins,$ let $\peM\za=(\pe\za)^\rL$. 

The set 
$\peim=\peim_\bI=
\bigcup_{\za\in\ins}\peM\za\in \rL$ 
will be the {\it forcing notion\/}. 

To define the order, we  put 
$\supp X=\za$ whenever $X\in\peM\za$. 
Now we set $X\leq Y$ 
(\ie\ $X$ is {\it stronger\/} than $Y$) 
iff $\za=\supp Y\sq \supp X$ and $X\res\za\sq Y$. 

\bre
\lam{bon}
We may note that the set $\bon=\ans\pu$ as in 
Remark~\ref{lao} belongs to $\peim$ and is the 
\dd\leq largest (\ie, the weakest) element 
of $\peim$.
\ere

Now let $G\sq\peim$ be a \dd\peim generic set 
(filter) over $\rL$. 
 
\bre
\lam{diez}
If $X\in\peM\za$ in $\rL$ then $X$ is not even  
a closed set in $\can\za$ in $\rL[G]$. 
However we can transform it to a perfect set in 
$\rL[G]$ by the closure operation. 
Indeed the 
topological closure $X^\#$ of such a set $X$ 
in $\can\za$ taken in $\rL[G]$ 
belongs to $\pe\za$ from the point of 
view of $\rL[G]$.
\ere

It easily follows from Lemma~\ref{clop} that there 
exists a unique array $\a[G]=\sis{\a_\i[G]}{\i\in\bI}$,  
all $\a_\i[G]$ being elements of $\dn$, 
such that $\a[G]\res\xi\in X^\#$ 
whenever $X\in G$ and $\supp X=\xi\in\ins$. 
Then $\rL[G]=\rL[\sis{\a_\i[G]}{\i\in\bI}]=\rL[\a[G]]$ 
is a \dd\peim{\it generic extension\/} of $\rL$.

\bte
[Theorems 24, 31 in \cite{fm97}]
\lam{card}
Every cardinal in\/ $\rL$ remains a cardinal in 
$\rL[G]$. 
Every\/ $\a_\i[G]$ is Sacks generic over the model\/ 
$\rL[\a[G]\rsd{\cls\i}]$.
\ete


We now present several lemmas on reals in   
\dd\peim generic models $\rL[G]$, established 
in \cite{fm97}. 
In the lemmas, we 
let $G\sq \dP $ be a set \dd\dP generic 
over $\rL$. 

\ble
[Lemma 22 in \cite{fm97}]
\lam{etaxi}
\label{ij}
Suppose that sets\/ $\eta,\,\xi\in\ins$ satisfy\/ 
$\kaz\j\in\eta\:\sus \i\in\xi\:(\j\cle\i)$. 
Then\/ $\a[G]\res\eta\in\rL[\a[G]\res\xi]$.
\ele

\ble
[Lemma 26 in \cite{fm97}]
\lam{xi.i}
Suppose that\/ $\bK\in\rL$ is an initial segment in\/ 
$\bI$, and\/ $\i\in\bI\setminus\bK$. 
Then\/ $\a_\i[G]\not\in\rL[\a[G]\res\bK]$.
\ele

\ble
[Corollary 27 in \cite{fm97}]
\lam{xixi}
If\/ $\i\not=\j$ then\/ $\a_\i[G]\not=\a_\j[G]$ 
and even\/ 
$\rL[\a_\i[G]]\not=\rL[\a_\j[G]]$.\qed
\ele 

\ble
[Lemma 29 in \cite{fm97}]
\lam{lm3}
If\/ $\bK\in \rL$ is an initial segment of\/ $\bI$,  
and\/ $r$ is a real in\/ $\rL[G]$, 
then either\/ $r\in\rL[\x\res\bK]$ or there 
is\/ $\i\not\in\bK$ such that\/ $\a_\i[G]\in\rL[r]$. 
\ele

We apply the lemmas in the proof of 
the next theorem.
Let $\lel$ denote the G\"odel wellordering on 
$\dn,$ so that $x\lel y$ iff $x\in\rL[y]$.  
Let $x\les y$ mean that $x\lel y$ but $y\not\lel x$, 
and $x\eql y$ mean that $x\lel y$ and $y\lel x$. 

\bte
\lam{ft}
Assume that\/ $\i\in\bI$ and\/ 
$r\in\rL[G]\cap\dn.$ 
Then
\ben
\renu
\itlb{ft1}
if\/ $\j\in\bI$ and\/ $\j\cle\i$ 
then\/ $\a_\j[G]\lel\a_\i[G]\;;$

\itlb{ft2}
if\/ $\j\in\bI$ and\/ $\j\not\cle\i$ 
then\/ $\a_\j[G]\not\lel\a_\i[G]\;;$

\itlb{ft3}
if\/ $r\lel\a_\i[G]$ then\/ $r\in\rL$ 
or\/ $r\eql\a_\j[G]$ 
for some\/ $\j\in\bI$, $\j\cle\i\;;$

\itlb{ft6}
if\/ $\i=\ang{\ga,s}\in\bI$, $e=0,1$, and\/ 
$\i\we e=\ang{\ga,s\we e}$ then\/ $\a_{\i\we e}[G]$ 
is a\/ {\rm\ubf true successor\/} of\/ $\a_\i[G]$ 
in the sense 
that\/ $\a_\i[G]\les\a_{\i\we e}[G]$ and any real 
$y\in\dn$ satisfies\/ 
$y\les\a_{\i\we e}[G]\imp y\lel\a_{\i}[G]\;;$

\itlb{ft7}
if\/ $\i=\ang{\ga,s}\in\bI$, and\/ 
$x\in\dn\cap\rL[G]$ 
is a true successor of\/ $\a_\i[G]$ 
in the sense of\/ \ref{ft6}, then there is  
$e=0$ or\/ $1$ such that\/ 
$x\eql\a_{\i\we e}[G]\;.$
\een
\ete
\bpf
\ref{ft1} 
Apply Lemma~\ref{etaxi} with $\eta=\ans\j$ 
and $\xi=\ans\i$.

\ref{ft2}
Apply Lemma~\ref{xi.i} with $\bK=[\cle\i]$.  

\ref{ft3}
If there are elements $\j\in\cI$, $\j\cle\i$, 
such that $\a_\j[G]\in\rL[r]$, then let $\j$ be the 
largest such one, and let $\xi=[\cle\j]$ 
(a finite initial segment of $\bI$).
Then, by Lemma~\ref{lm3}, 
either\/ $r\in\rL[\a[G]\res\xi]$, or there 
is\/ $\i'\not\in\xi$ such that\/ 
$\a_{\i'}[G]\in\rL[r]$. 

In the ``either'' case, we have $r\in\rL[\a_\j[G]]$ 
by \ref{ft1}, so that  $\rL[r]=\rL[\a_\j[G]]$ by the 
choice of $\j$. 
In the ``or'' case we have 
$\a_{\i'}[G]\in\rL[a_\i[G]]$, 
hence $\i'\cle\i$ by \ref{ft2}. 
But this contradicts the choice of $\j$ and $\i'$. 

Finally if there is no $\j\in\cI$, $\j\cle\i$, 
such that $\a_\j[G]\in\rL[r]$, then the same argument 
with $\xi=\pu$ gives $r\in\rL$.

\ref{ft6}
The relation $\a_\j[G]\les\a_{\i\we e}[G]$ is 
implied by Lemmas \ref{etaxi} and \ref{xi.i}.
If now $y\les\a_{\i\we e}[G]$ then 
$y\in\rL$ or 
$y\eql \a_\j[G]$ for some $\j\cle\i\we e$ by 
\ref{ft3}, and in the latter case in fact 
$\j\cls\i\we e$, hence $\j\cle\i$, and 
then $y\lel\a_{\i}[G]$.

\ref{ft7}
By \ref{ft6}, it suffices to prove that 
$x\lel\a_{\i\we 0}[G]$ or $x\lel\a_{\i\we 1}[G]$.
Assume that $x\not\lel\a_{\i\we 0}[G]$. 
Then by 
Lemma~\ref{lm3} there is an element $\j\in\bI$ 
such that $\j\not\cle\i\we0$ and 
$\a_\io[G]\lel x$. 
If $\a_\j[G]\les x$ strictly then 
$\a_\j[G]\lel\a_\i[G]$ by the true successor 
property, hence $\io\cle\i$, contrary to 
$\io\not\cle\i\we0$, see above. 
Therefore in fact $\a_\io[G]\eql x$. 
Then we must have 
$\io=\i\we0$ or $\io=\i\we1$ as $x$ is a true 
successor, but then $\io=\i\we1$, as 
$x\not\lel\a_{\i\we 0}[G]$ was assumed, 
and we are done. 
\epf

\parf{The subextension}
\las{subex}

Following the arguments above, assume that 
$G\sq \dP $ is a set \dd\dP generic 
over $\rL$,   
and consider the set 
$\jg\in\rL[G]$ of all 
elements $\i\in\bI$ such that either 
$\i=\ang{\ga,0^m}$, where $\ga<\omi$ and $m<\om$, 
or $\i=\ang{\ga,0^m\we1}$, where $\ga<\omi$ and 
$m<\om$, $\a_\ga[G](m)=1$.
Following \eqref{eq4}, we define 
\busq
{eq5}
{\mg\;=\;\pws\om\:\cap\:
\bigcup_{\i_1,\dots,\i_n\in\jg} 
\rL[a_{\i_1}[G],\dots,a_{\i_n}[G]],
}  

\ble
\lam{77}
If\/ $\i\nin\jg$ then\/ $\a_\i[G]\nin\mg$.
\ele
\bpf
This is not immediately a case of Lemma~\ref{xi.i} 
because $\jg\nin\rL$. 
However the set $\bK=\ens{\j\in\bI}{\i\not\cle\j}$ 
belongs to $\rL$ and satisfies  
$\jg\sq\bK\sq\bI$.
We have $\i\nin\bK$, and hence 
$\a_\i[G]\nin\rL[\a[G]\res\bK]$ by Lemma~\ref{xi.i}. 
On the other hand, we easily check  
$X\sq \rL[\a[G]\res\bK]$, 
and we are done.
\epf
 
We are going to prove that 
$\stk\om\mg$ is a model of $\pads+\CA(\us12)$, 
but the full $\CA$ fails in $\stk\om\mg$.\vom

{\ubf Part 1:} 
$\stk\om\mg$ is a model of all axioms of $\pad$ 
except for $\CA$, trivial.\vom

{\ubf Part 2:} 
$\stk\om\mg$ is a model of $\CA(\us12)$
(with parameters). 
This is also easy by the Shoenfield absoluteness 
theorem.\vom

{\ubf Part 3:} 
$\stk\om\mg$ fails to satisfy the full $\CA$. 
Here we need some work.
Let $\ga<\omil$, so that both $\ga$ and each 
pair $\ang{\ga,s}$, $s\in\bse,$ belong to $\bI$ 
by \eqref{eq2} in Section~\ref{gsi}, in particular 
$\i_0=\ang{\ga,\La}\in\bI$, 
where $\La$ is the empty tuple.
In addition $\ga$ (as an element of $\bI$) does 
not belong to $\jg$. 
Our plan is to prove that $\a_\ga[G]\nin\mg$ but 
$\a_\ga[G]$ is definable in $\stk\om\mg$.\vom

{\ubf Subpart 3.1:} $\a_\ga[G]\nin\mg$ by 
Lemma~\ref{77} just 
because $\ga\nin\jg$.\vom

{\ubf Subpart 3.2:} 
$\a_\ga[G]$ is definable in $\stk\om\mg$ with 
$\a_\io[G]$ as a parameter, 
where $\io=\ang{\ga,\La}\in\jg$. 
Namely we claim that for any $m<\om$: 
\pagebreak[0] 
\busq
{eq6}
{\a_\ga[G](m)=1 
\quad\text{iff}\quad
\begin{minipage}[t]{0.63\textwidth}
there is an array of reals 
$b_0,b_1,\dots,b_m,b_{m+1}$ and $b'_{m+1}$ 
in $\dn$ such that $b_0=\a_\io$, 
each $b_{k+1}$ is a true successor of $b_k$ 
($k\le m$), $b'_{m+1}$ is a true successor 
of $b_m$ as well, and $b'_{m+1}\not\eql b_{m+1}$.
\end{minipage}
}

The formula in the right-hand side of 
\eqref{eq6} is based on the G\"odel canonical 
$\is12$ formula for $\lel$, which is absolute 
for $\mg$ by the definition of $\mg$. 
Therefore \eqref{eq6} implies that  
$\a_\ga[G]$ is definable in $\stk\om\mg$ with 
$\a_\io[G]$ as a parameter. 
Thus it remains to establish \eqref{eq6}. 

Direction $\imp$. 
Assume that $\a_\ga[G](m)=1$. 
Then $\jg$ contains the elements 
$\i_k=\ang{\ga,0^k}$, $k\le m+1$, along with 
an element $\i'_{m+1}=\ang{\ga,0^m\we 1}$.
Therefore the reals 
$b_k=\a_{\i_k}[G]$, $k\le m+1$, 
and $b'_{m+1}=\a_{\i'_{m+1}}[G]$ belong to $\mg$. 
Now Theorem~\ref{ft}\ref{ft6},\ref{ft2} implies 
that the reals $b_k$ and $b'_{m+1}$ satisfy 
the right-hand side of \eqref{eq6}, as required.

Direction $\mpi$. 
Assume that the reals $b_k$, $k\le m+1$, 
and $b'_{m+1}$ satisfy 
the right-hand side of \eqref{eq6}. 
By Theorem~\ref{ft}\ref{ft7}, there is an 
array of bits $e_1,\dots,e_m,e_{m+1}$ and 
$e'_{m+1}$ such that 
$b_k=\a_{\i_k}[G]$ for all $k\le m+1$ and 
$b'_{m+1}=\a_{\i'_{m+1}}[G]$, where 
$\i_k=\ang{\ga,\ang{e_1,\dots,e_k}}$ and 
$\i'_{m+1}=\ang{\ga,\ang{e_1,\dots,e_m,e'_{m+1}}}$. 

However we must have 
$\i_k\in\jg$ for all $k\le m+1$, and 
$\i'_{m+1}\in\jg$, by Lemma~\ref{77}, 
since the reals $b_k$ 
and $b'_{m+1}$ belong to $\mg$. 
Then obviously $e_1=\dots=e_m=0$ while 
$e_{m+1}=0$ and $e'_{m+1}=1$ or vice versa 
$e_{m+1}=1$ and $e'_{m+1}=0$. 
In other words, the elements 
$\ang{\ga,0^{m+1}}$ and $\ang{\ga,0^{m}\we1}$
belong to $\jg$. 
This implies $\a_\ga[G](m)=1$.\vom

{\ubf Part 4:}
$\stk\om\mg$ satisfies the \paf\ schema $\CAs$. 
This is rather similar to the verification of 
$\CAs$ in $\stk\om\xg$ in 
Section~\ref{coh}. 

Assume that $\Phi(k)$ is a 
\paf\ $\lpad$ formula with $k$ the only 
free variable. 
Consider the set 
$y=\ens{k<\om}{\stk\om\mg\mo\Phi(k)}$; 
then $y\in\rL[G]$, $y\sq\om$. 
We claim that $y$ even belongs to $\rL$, 
and hence to $\mg$. 

Let $\ofo$ be the forcing relation associated 
with $\peim$, over $\rL$ as the ground model. 
Thus if $X\in\peim$ and $k<\om$ then 
$X\ofo\Phi(k)$ iff $\Phi(k)$ holds in any 
\dd\dP generic extension $\rL[H]$ of $\rL$ 
such that $X\in H$.\snos{ 
See Kunen~\cite{kun} on forcing, 
especially Section IV.6 there 
on the ``forcing over the universe'' approach.}
Let $\uG$ be a canonical \dd\peim name for $G$.
We assert that 
\busq
{eq7}
{
y=\ens{k<\om}
{\bon\ofo\text{``$\stk\om\mug\mo \Phi(k)$''}}.
}
(See Remark~\ref{bon} on $\bon$.)

In the nontrivial direction, assume that $k\in y$. 
Then by the forcing theorem there is a condition 
$X\in G$ forcing $\stk\om\mug\mo \Phi(k)$. 
We claim that then $\bon$ forces the same as well. 

To prove this reduction, we define, 
{\ubf still in $\rL$}, 
the set $\Perm\in\rL$ that consists of all 
bijections $\pi:\omi\onto\omi$ such that 
$\pi=\pi\obr$ and 
\rit{the domain of nontriviality} 
$\abc \pi=\ens{\al}{\pi(\al)\ne\al}$ is 
at most countable, \ie, bounded in $\omi$. 
Any $\pi\in\Perm$ acts on:
\bit
\item[$-$]
elements $\i=\ga$ or $\i=\ang{\ga,s}$ of 
$\bI$, by $\pi\i=\pi(\ga)$, resp.\ 
$\i=\ang{\pi(\ga),s}$;

\item[$-$]
maps $g$ with $\dom g\sq\bI$, by 
$\dom(\pi g)=\ima\pi{\dom g}$ and 
$(\pi g)(\pi(\al))=g(\al)$ for all 
$\al\in\dom g$;

\item[$-$]
thus if $\xi\sq\bI$ and 
$x\in\can\xi$ then $\pi x\in\can{\ima\pi\xi}$ 
and $(\pi x)(\pi(\al))=x(\al)$;

\item[$-$]
sets $X\in\pe\xi$, $\xi\in\cpo$, by 
$\pi X=\ens{\pi x}{x\in X}\in\pe{\ima\pi\xi}$.
\eit

We return to the nontrivial direction $\imp$ 
of \eqref{eq7}, where we have to prove that 
the condition $\bon$ forces  
``$\stk\om\mug\mo \Phi(k)$''. 
Let this be not the case. 

Then there is a condition $Y\in\peim$ 
which forces ``$\stk\om\mug\mo \neg\,\Phi(k)$''. 
There is a permutation $\pi\in\Perm$ satisfying 
$\supp{Z}\cap\supp X=\pu$, where $Z=\pi Y\in\peim$. 
We claim that $Z$ forces 
``$\stk\om\mug\mo \neg\,\Phi(k)$''. 
Indeed assume that $H\sq\peim$ is a set 
\dd\peim generic 
over $\rL$, and $Z\in H$. 
We have to prove that 
$\stk\om{\mh}\mo\neg\,\Phi(k)$.
The set $K=\ens{\pi{Z'}}{Z'\in H}$ is \dd\bP generic 
over $\rL$ along with $H$ since $\pi\in\rL$. 
Moreover $K$ contains $Y$. 
It follows that 
$\stk\om{\mk}\mo\neg\,\Phi(k)$
by the forcing theorem and the choice of $Y$. 

However the array $\a[K]$ is equal to 
the permutation of the array $\a[H]$ by $\pi$. 
It follows that $\mh=\mk$, and hence 
$\stk\om{\mh}\mo\neg\,\Phi(k)$, as required. 
Thus indeed $Z$ forces 
``$\stk\om\mug\mo \neg\,\Phi(k)$''. 

Recall that $X$ forces ``$\stk\om\mug\mo \Phi(k)$''. 
On the other hand, $X,Z$ are compatible in $\peim$ 
because $\supp{Z}\cap\supp X=\pu$. 
This is a contradiction. 

We conclude that $\bon$ forces 
``$\stk\om\mug\mo \Phi(k)$'', 
and this completes the proof of \eqref{eq7}. 
But it is known that the forcing relation $\ofo$ 
is expressible in $\rL$, the ground model. 
Therefore it follows from \eqref{eq7} 
that $y\in\rL$, hence $y\in \mg$, as required.

\parf{Remarks and questions}
\las{rq}

Here we present three questions related to possible 
extensions of Theorem~\ref{mt2}.

\bvo
\lam{81}
Is the \paf\ countable choice schema $\AC^\ast$ 
in the language $\lpad$ 
true in the models $\stk\om\mg$ defined in 
Section~\ref{subex}\:?
\evo
  
\bvo
\lam{82}
Can we sharpen the result of Theorem~\ref{mt2} 
by specifying that $\CA(\us13)$ is violated? 
The combination $\CA(\us12)$ plus $\neg\,\CA(\us13)$
would be optimal. 
The counterexample to $\CA$ defined in 
Section~\ref{subex} (Part 3) definitely is 
more complex than $\us13$. 
\evo
  
\bvo
\lam{83}
As a generalization of the above, prove that, 
for any $n\ge2$,  
$\pads + \CA(\us1n)$ does not imply  
$\CA(\us1{n+1})$. 
In this case, we'll be able to 
conclude that the full schema $\CA$ 
is not finitely axiomatizable over $\pads$.
Compare to Problem 9 in \cite[\S\,11]{aptm}.
\evo

\back
The authors are thankful to Ali Enayat for his 
enlightening comments that made it possible to 
accomplish this research. 
\vyk{
The authors are thankful to 
the anonymous referee for their 
comments and suggestions, which significantly 
contributed to improving the quality of the publication.
}
\eack


\renek{\refname}
{{References}\addcontentsline{toc}{section}{References}}
\small

\bibliographystyle{plain} 

\bibliography{69.bib,69kl.bib}


\end{document}